\documentclass[11pt]{article}
\usepackage{cite}
\usepackage{mathrsfs}
\usepackage{amsfonts}
\usepackage{amsmath}
\usepackage{amsfonts,amssymb,color}
\usepackage{dsfont}
\usepackage{curves}
\usepackage{mathrsfs}
\usepackage{pifont}
\usepackage{amssymb}
\usepackage{latexsym,amsmath,amssymb,amsfonts,epsfig,graphicx,cite,psfrag}
\usepackage{eepic,color,colordvi,amscd}

\newtheorem{theorem}{Theorem}[section]

\newtheorem{lemma}[theorem]{Lemma}
\newtheorem{corollary}[theorem]{Corollary}

\newcommand{\qed}{\hfill\rule{0.5em}{0.809em}}

\numberwithin{figure}{section}
\numberwithin{equation}{section}

\def\emptyset{\mbox{{\rm \O}}}

\def\qed{\hfill \rule{4pt}{7pt}}
\def\pf{\noindent {\it Proof. }}

\date{}
\begin{document}
\title{Thickness and Outerthickness for Embedded Graphs}
 \author{\small   Baogang  Xu\footnote{Research partially supported by NSFC projects 11331003 and 11571180, and by a project funded by PAPD of Jiangsu Higher Education Institutions.}\\
\small Institute of Mathematics, School of Mathematical Sciences\\
\small Nanjing Normal University, 1 Wenyuan Road,  Nanjing, 210023,  China\\
\small baogxu@njnu.edu.cn\\
\small Xiaoya Zha\footnote{Research  supported by National Security Agency grant H98230-13-1-0216.}\\
\small Department of Mathematical Sciences\\
\small Middle Tennessee State University, Murfreesboro, TN 37132, U.~S.~A.\\
\small xzha@mtsu.edu
}

\maketitle

\begin{abstract}
We consider the thickness $\theta(G)$ and outerthickness $\theta_o(G)$ of a graph $G$ in terms of its orientable and nonorientable genus. Dean and Hutchinson provided upper bounds for thickness of graphs in terms of their orientable genus. More recently, Concalves proved that the outerthickness of any planar graph is at most $2$. In this paper, we apply the method of deleting  spanning disks of embeddings to approximate the thickness and outerthickness of graphs. We first obtain  better upper bounds for thickness.
We then use a similar approach to provide upper bounds for outerthickness of graphs in terms of their orientable and nonorientable genera. Finally we show that the outerthickness of the torus (the maximum outerthickness of all toroidal graphs) is $3$. We also show that all graphs embeddable in the double torus have thickness at most $3$ and outerthickness at most $5$.
\begin{flushleft}
{\em Key words and phrases:} surface; planar graph; outerplanar graph; decomposition.\\
{\em AMS 2000 Subject Classifications:}  05C10, 05C35\\
\end{flushleft}

\end{abstract}
\newpage
\section{Introduction and Terminology}

An {\em outerplanar graph} is a planar graph that can be embedded in the plane without crossing edges,
in such a way that all the vertices are incident with the same face.
The {\em thickness} of a graph $G$, denoted by $\theta(G)$  (first defined by Tutte \cite{Tut63}),  is the minimum number of planar subgraphs whose union is $G$.
Similarly, the {\em outerthickness} $\theta_o(G)$ is obtained  where the planar subgraphs are replaced by outerplanar subgraphs in the previous definition. If $\Sigma$ is a surface, define $\theta(\Sigma) = \max \{\theta(G)$ :  $G$ is embeddable in $\Sigma \}$, where the maximum is taking over all graphs embeddable in $\Sigma$. Define $\theta_o(\Sigma) $ analogously.

Much work has been done in partitioning the edges of graphs such that each subset induces a subgraph of a certain type. A well-known result by Nash-Williams \cite{Nash} gives a necessary and sufficient condition for a graph to admit an edge-partition into a fixed number of forests. His results imply that any planar graph can be edge-partitioned into three forests, and any outerplanar graph into two forests.
A lot of research has been devoted  to partitioning the edges into planar graphs (to determine the thickness of graphs) and outerplanar graphs (to determine the outerthickness of graphs).
The thickness of several special classes of graphs have been determined, including the complete graphs $K_n$ \cite{AG76, JV76} (see (\ref{thick-complete})), the complete bipartite graphs $K_{m, n}$  \cite{BHM64} (except possibly if $m$ and $n$ are both odd, or $m\le n$ and $n$ takes some special values), and  the hypercube $Q_n$ \cite{MK1967}. See the survey paper \cite{MOS98} for more results on thickness of graphs. Guy and Nowakowski \cite{GN901, GN902} determined  the outerthickness of complete graphs (see (\ref{outer-complete})),  the hypercube and some complete bipartite graphs. Here,
\begin{equation}\label{thick-complete}
\theta(K_{n})= \left\{\begin{array}{ll}
\lfloor{n+7\over 6}\rfloor, & \; \mbox{ if } n\neq 9, 10,\\
3, & \; \mbox{ otherwise, }\end{array}\right.
\end{equation}
and
\begin{equation}\label{outer-complete}
\theta_o(K_n)=\left\{\begin{array}{ll}
\lceil{n+1\over 4}\rceil, & \; \mbox{ if } n\neq 7\\
3, & \; \ n=7 \end{array}\right.
\end{equation}

It is known that thickness problem is ${\cal N}{\cal P}$-hard \cite{Man83}, and there are no  other classes of graphs for which the thickness have been found. For more general classes of graphs, the attention has been focused on  finding upper bounds of thickness and outerthickness.
J\"unger et al \cite{JMOS98} have shown that a graph has thickness at most $2$ if it contains no $K_5$-minor.
Asano \cite{As87} proved that, if a graph $G$ is triangle free and has orientable genus $\gamma$,  then $\theta(G) \le \gamma(G) + 1$. He also showed that all toroidal graphs have thickness at most $2$.
Dean and Hutchinson \cite{DH88} strengthened Asano's result by proving that $\theta(G) \le 6 + \sqrt{2 \gamma (G) -1}$.

In 1971, Chartrand, Geller and Hedetniemi \cite{CGH71} conjectured that every planar graph has an edge partition into two outerplanar graphs. Ding, Oporowski, Sanders  and Vertigan \cite{DOSV00} proved that every planar graph has an edge partition  into two outerplanar graphs and a {\em vee-forest}, where a vee-forest is the disjoint union of a number of $K_2$'s and $K_{1, 2}$'s. They also showed that every graph with nonnegative Euler characteristic has an edge partition into two graphs of tree-width at most three.
Recently Gongalves \cite{Go05, Go08} confirmed the Chartrand, Geller and Hedetniemi's conjecture by showing that the edge set of a planar graph can be partitioned into two outerplanar graphs.
An interesting result by Kedlaya \cite{Ke96} shows that some planar graphs cannot be edge-partitioned into two outerplanar subgraphs such that one of them is outerplanarly embedded.

In this paper, we first provide some technical results in Section 2.
We introduce the technique of deleting maximal spanning disks for embeddings of graphs.
We also introduce  essential edges with respect to spanning disks of embeddings, and the corresponding noncontractible  nonhomotopic loop system of surfaces.
Applying these techniques we provide results on thickness and outerthickness to
(i) improve Dean and Hutchinson's upper bounds for graphs in terms of their orientable  and nonorientable genus (Section 3), (ii) obtain upper bounds for outerthickness of graphs in terms of their orientable and nonorientable genus  (Section 3), (iii) show that  outerthickness of the torus is $3$ (Section 4),
  and (iv) improve Asono's result by dropping his triangle free condition, and show that all graphs embeddable in double torus have thickness at most $3$ and outerthickness at most $5$ (Section 4),
      and all graphs embedded in triple torus have thickness at most $4$ (Section 4) .

\section{Technical results}

We prove some technical and structural results in this section.
Since adding multiple edges to a graph $G$ does not increase the thickness/outerthickness of $G$,
and the thickness/outerthickness of a graph is equal to the maximum thickness/outerthickness of its blocks,
we may assume that graphs are simple and $2$-connected. Let $S_g$ be the orientable surface with genus $g$ ($g\ge 0$, the sphere with $g$ handles) and $N_k$ be the nonorientable surface with nonorientable genus $k$
($k\ge 1$, the sphere with $k$ crosscaps).
Suppose $C$ is a cycle of a graph embedded in surface $\Sigma$, and $x$ and $y$ are two vertices on $C$.
We assign a direction to $C$ and define $xCy$ to be the open path from $x$ to $y$ in this direction.
The following is obvious.

\begin{lemma}\label{lem-2-1}
If $G$ is a subgraph of $H$, then $\theta(G) \le \theta(H)$ and $\theta_o(G) \le \theta_o(H)$.
\end{lemma}

In order to study thickness/outerthickness of graphs embedded in surfaces,
we will apply Lemma~\ref{lem-2-1} by adding edges to $G$ to obtain a spanning supergraph $H$ of $G$ then study the thickness/outerthickness of $H$. In this way we may obtain a better structure of embeddings.
Note that Lemma~\ref{lem-2-1} may not be true if the subgraph relation is replaced by the minor or subdivision relations.

Let $G$ be a  graph and $\Psi(G)$ be an embedding of $G$ in a surface $\Sigma$.
A subembedding $\Psi^s$ is {\em spanning} if it contains all vertices of $G$.
A spanning subembedding is {\em contractible} if it does not contain any noncontractible cycle of $\Psi(G)$.
In particular a contractible spanning subembedding is a {\em spanning disk} if it is homeomorphic to a closed disk, in which case the boundary of this spanning subembedding is a contractible cycle of $\Psi(G)$.
For any embedding,  a spanning tree is always a contractible spanning subembedding.
However, an embedding may not contain a spanning disk.
An example is the unique embedding of the Heawood graph in the torus which is the dual embedding of $K_7$.
It contains no spanning disk even though the embedding has face width three (or equivalently, polyhedral embedding, or wheel-neighborhood embedding). An edge $e$ is {\em essential}, with respect to a contractible spanning subembedding $\Psi^s$ if $e \cup \Psi^s$ contains a noncontractible cycle.
Note that if $e$ is an essential edge then $e$ is contained in every noncontractible cycle of $e \cup \Psi^s$.
An essential edge becomes a noncontractible loop if we contract $\Psi^s$ to a single point.

\begin{lemma}\label{lem-2-2}
 Let $G$ be a simple graph and $\Psi(G)$ be an orientable genus embedding
or a minimal surface  embedding (with maximum Euler Characteristic) of $G$ in  $\Sigma$.
   Then $G$ has a simple spanning supergraph $H$ embedded in $\Sigma$ with embedding $\Psi(H)$ such that
 $\Psi(H)$  is an orientable genus embedding or a minimal surface  embedding (which must be cellular) and  contains a spanning disk.
\end{lemma}

Lemma~\ref{lem-2-2} is not true if the embedding is not an orientable genus embedding or a minimal surface embedding. For example, if $G$ is a  complete graph embedded in its maximal surface $\Sigma$ then there does not exist an $H$ embedded in $\Sigma$ with $G$ being a spanning subgraph of  $H$ such that the embedding contains a spanning disk.

\medskip

\noindent{\bf Proof of Lemma}~\ref{lem-2-2}.
Let $\Psi(G)$ be a minimal surface embedding of $G$ in $\Sigma$.
We start with a spanning tree $T$ of $G$, and add more faces to $T$ such that the resulting subembedding $R$ is maximal and contractible.
We then add new edges to $R$ one by one such that the resulting graph is a supergraph $G^+$ of $G$ with the same vertex set, and has a spanning region $R^+$ which is contractible.
We do not add multiple edges to $G$, and thus $G^+$ remains to be a simple graph.
We assume that $R^+$ is constructed such that it has minimum number of essential edges.
We now claim that $R^+$ is a spanning disk of $G^+$.

Suppose not. Since $R^+$ is contractible, there exists a vertex $v$ which is a cut vertex of $R^+$. Note that $v$ is not a cut vertex of $G^+$ as $G^+$ is $2$-connected.
Therefore there is a face $f$ of $\Psi$ such that $v$ appears at least twice on the boundary of $f$,
and there exists a noncontractible simple closed curve $\Gamma$ of $\Sigma$ such that
$\Gamma$ intersects $G^+$ only at $v$.
Let $x$ and $y$ be the two vertices incident to $v$ which appear on the facial walk of $f$ in the order of $xvy$.
If $xy$ is not an existing edge then we add $xy$ along with the facial walk $xvy$ to obtain a new supergraph.
Now add the face between $xy$ and $xvy$ to enlarge $R^+$ to obtain a new contractible region of the new supergraph with fewer cut vertices of the spanning region.
If $xy$ is already an existing edge, then $xy \cup xvy$ is a cycle of $G^+$ intersecting $\Gamma$ only at $v$.
Therefore $xy \cup xvy$ is a noncontractible cycle and thus $xy$ is an essential edge of $R^+$.
We re-embed $xy$ along $xvy$ to obtain a new embedding with one less essential edge, contradicting the assumption that $R^+$ has minimum number of essential edges.
The new embedding is a cellular embedding if $\Psi(G)$ is an orientable genus embedding or a minimal embedding. For otherwise, if the embedding is not cellular embedding after re-embeding of $xy$, then the original face incident to $xy$ must have contained a noncontractible simple closed curve $\Gamma_1$ which is homotopically disjoint from $G^+$.
Cutting the surface $\Sigma$ along $\Gamma _1$ and capping off the boundary would then result in an embedding of $G$ in a surface with higher Euler characteristic, a contradiction.
Therefore the resulting embedding is cellular and  $R^+$ is a spanning disk of $G^+$.
This proves  Lemma~\ref{lem-2-2}.  \qed

\medskip

Let $\Psi(G)$ be an embedding  of a graph $G$  in  a  surface $\Sigma$.
We now allow $G$ to have multiple edges,
but if $G$ has multiple edges then any two multiple edges form a noncontractible cycle
(this is to prevent two multiple edges from forming a face of size $2$).
Suppose $\Psi(G)$ has a spanning disk $D$.
Denote the subgraph embedded in $D$ by $D(G)$ (including all edges on the boundary of $D$).
Then $\Psi(G) \setminus D(G)$ consists of  essential edges only, which we call {\em subembedding of essential edges}. Denote this subembedding by $G^e$. We also use $G^e$ to represent the subgraph consisting of all essential edges.

 Let $\Sigma$ be a surface with Euler characteristic $\chi (\Sigma) $.
       Assume   $x$ is a point on  $\Sigma$.
       Let ${\cal L} =\{ l_i\ : \ i=1,2,...,t\}$ be a collection of
        noncontractible   loops with base point $x$  such that
       $l_i$ and $l_j$ only intersect at $x$ and  are not homotopic to each other for  $1 \le i < j \le t$.
       Call ${\cal L}$ a {\em nonhomotopic loop system}. A nonhomotopic loop system ${\cal L}$ is {\em maximal} if adding any noncontractible loop  $l$ with $x$ as the base point to ${\cal L}$ then $l$ will either be  homotopic to some loop $l_i$ in  ${\cal L}$ or intersect  $l_i$ at $x$ and also some other points.
       Let $\rho(\Sigma) = \max \{ |{\cal L}| \} $, where ${\cal L}$ is a maximal nonhomotopic loop system of $\Sigma$,  $|{\cal L}|$ is the number of loops in ${\cal L}$, and the the maximality is taken over all such systems of $\Sigma$.

A maximal nonhomotopic loop system of $\Sigma$ is closely related to the structure of the subembedding of essential edges $G^e$ since all essential edges become noncontractible loops when we contract the spanning disk $D$ to a single point. Not many results and references can be found for the structure of a maximal nunhomotopic loop system of surfaces and its maximum number. However, we will determine $\rho(S_1)$ in Section 4 which helps us to determine $\theta_o(S_1)$.

\section{Upper bounds for thickness and outerthickness for general surfaces}

In this section we study the upper bounds for thickness and outerthickness for general surfaces.
We will improve both Asano, and Dean and Hutchinson's results on thickness of graphs in terms of their genus.
We also obtain similar upper bounds for outerthickness of graphs  in terms of their genus.
As mentioned in the introduction, Asano \cite{As87} proved that, if $G$ has orientable genus $\gamma$, then  $\theta(G) \le \gamma(G) +1$ if $G$ has no triangle.
This result is certainly weaker than the later result by Dean and Hutchinson \cite{DH88} who showed that $\theta(G) \le 6 + \sqrt{2\gamma -2}$.
However the $\gamma(G) +1$ bound works better for surface with lower genus.
For example, for triangle-free graphs, Asano's bound for the double torus is $3$, but Dean and Hutchinson's bound for the double torus is $7$. Our first result in this section is to show that one may drop the triangle free assumption for Asano's bound.

\begin{theorem}\label{thm-3-1}
If $G$ is embedded in its orientable genus surface $S_{\gamma}$  then $\theta(G) \le \gamma (G) +1$.
\end{theorem}
\pf We prove the theorem by induction on $\gamma(G)$.
If $G$ is planar then $\gamma(G) =0$ and $\theta(G) =1$, so the theorem is true for $\gamma(G)=0$.
Assume that the theorem is true if $\gamma(G) \le n-1$.
Now assume that $\gamma(G) = n$ and let $\Psi(G)$ be an embedding of $G$ in $S_{n}$.
By Lemmas~\ref{lem-2-1} and \ref{lem-2-2}, we may assume that $\Psi(G)$ has a spanning disk $D$.
Let $G^e$ be the subembedding of essential edges, and $v$ be a vertex of $G^e$.
  Let {\em St($v$)} be the subgraph of $G^e$ consisting of  all essential edges incident to $v$.
  Then St($v$) is a star.
  Let $vu$ be the first essential edge in the clockwise rotation of all essential edges incident to $v$ (here the first edge means the edge immediately following the edge which is incident to $v$  on the boundary of $D$).
  Let $G_1 = D \cup {\it St}(u)$.
  Then $G_1$  is planar because, if necessary, we can connect all essential edges incident to $u$ to the boundary of $D$ by permuting the clockwise order of
  essential edges incident to $u$ inherited from $\Psi(G)$.

  Let $L$ be a noncontractilbe cycle of $G$ containing $vu$ which is contained in $vu \cup D$.
  All noncontractible cycles containing $vu$ and are contained in $vu \cup D$ are homotopic.
  Let $\Gamma$ be a noncontractible simple closed curve of $S_{n}$ that is homotopic to $L$.
  The curve $\Gamma$ can be chosen such that $\Gamma \cap (G^e \backslash {\it St}(v) )= \emptyset$
  because $vu$ is the first edge of ${\it St}(v)$. Then $L$ is orientation preserving,
   and  all essential edges incident to $v$ are on the same side of $vu$.

  If $\Gamma$ is a noncontractible separating simple closed curve,
  cut $S_{n}$ along $\Gamma$ and cap off the two holes to obtain two orientable surfaces of genus $a$ and $b$, respectively, with $1 \le a, b \le n -1$ and $a+b =n$.
  Let $H_a$ and $H_b$ be the subgraphs of $G^e$ that are embedded in $S_a$ and $S_b$ respectively.
  By induction hypothesis, $\max \{ \theta (H_a) , \theta (H_b) \} \le (n-1) +1 =n$.
  Since $H_a\cap H_b = \emptyset$, we have $\theta(G^e) = \max \{\theta(H_a), \theta(H_b)\} \le (n-1)  +1=n$.

  Therefore $\theta(G) \le n +1 = \gamma(G) +1$, and the theorem is true in this case.

  If $\Gamma$ is a noncontractible nonseparating simple closed curve,
   cut $S_{n}$ along $\Gamma$ and cap off the two boundary  components with disks to obtain a surface which is $S_{n -1}$.
  Therefore $G^e \backslash {\it St}(u)$ is a subgraph of $G^e$ which is embedded in the surface $S_{n -1}$.
  By induction hypothesis, $\theta(G^e \backslash {\it St}(u)) \le (n-1)+1 =n$
  and thus $\theta(G) \le n+1 = \gamma(G) +1$.
  This completes the proof of Theorem~\ref{thm-3-1}.   \qed

\medskip

By induction on the genus and starting with $\theta (S_0) =1$ we have

\begin{corollary}\label{coro-3-2}
(i) $\theta(S_k) \le k+1, k\ge 1$, (ii) $\theta(S_2 ) \le 3$, and (iii) $\theta(S_3) \le 4$.
\end{corollary}

The upper bounds for thickness in Corollary~\ref{coro-3-2} are better for $k=2$ and 3 than the bounds obtained by Theorem~\ref{thm-3-4} later in this section,
which show that $\theta(S_2) \le 4$ and $\theta(S_3) \le 5$.

With the same arguments using induction on genus and $\theta_o(S_0) =2$ we have

\begin{corollary}\label{coro-3-3}
For $k \ge 1, \theta _o (S_k) \le \theta_o (S_{k-1}) +2$ and $\theta _o (S_k) \le 2k+2$.
\end{corollary}


A  result similar to Theorem~\ref{thm-3-1} can be obtained in terms of the Euler Characteristic.
However, the statement may be a little more complicated because when one reduces surface to a surface with higher Euler characteristic, the  Euler characteristic may increase either by $1$ or $2$,
depending on whether the cutting noncontractible simple closed curve is
orientation   reversing or orientation preserving.
We choose not to present such a theorem since later theorems provide better bounds.

\medskip

We now try to improve Dean and Hutchinson's upper bounds for thickness in terms of the genus.

The following statement is a part of the proof of Theorem 3 in \cite{DH88}.
We present here a new proof.

\begin{lemma}\label{lem-3-1}
Let $G$ be a graph, $d$ be a positive integer, and let $H$ be a graph obtained from $G$ by iteratively removing vertices of degree at most $d$.
Then, the edges in $E(G)\setminus E(H)$ can be decomposed into $d$ forests $F_1, \ldots, F_d$ (possibly empty) such that, for each $i$, each component of $F_i$ has at most one common vertex with $H$.
\end{lemma}
\pf The lemma is trivially true if $H=G-v_1$ for some vertex $v_1$, as we may choose each edge incident with $v_1$ as a forest.  Suppose that the conclusion holds while $H'$ is obtained from $G$ by removing at most $k-1$ vertices,  $k\ge 2$, and let $F'_1, \ldots, F'_d$ be a decomposition of $E(G)\backslash E(H')$ into forests as required.

Now, let $H=H'-v_k$, where $d_{H'}(v_k)\le d$.
Without loss of generality, we suppose that $d_{H'}(v_k)= d$, and let $v_ku_1, v_ku_2, \ldots, v_ku_d$ be the edges incident with $v_k$ in $H'$. For each $i\in \{1,2,\ldots, d\}$, let $F_i=F'_i+v_ku_i$.
Let $C$ be a component of $F_i$.
Then, either $V(C)\cap V(H)=\{u_i\}$ if $C$ is a subtree obtained from the components of $F'_i$ which contains either $u_i$ or $v_k$ (but not both by the
induction hypothesis) by adding $v_ku_i$, or $V(C)\cap V(H)=V(C)\cap V(H')$ if $C$ is a component of $F'_i$ which contains neither $u_i$ nor $v_k$.
In both cases, we see that $|V(C)\cap V(H)|\le 1$. This proves the lemma. \qed

\medskip

The proofs of the following two theorems are similar to those of \cite{DH88}, but employ a special embedding (i.e., an embedding with a spanning disk) of  graphs on surfaces.

\begin{theorem}\label{thm-3-4}
Let $G$ be a connected  graph embedded in the surface $S_g$ ($g\ge 1$).
  Then, $\theta(G)\le 3 + \sqrt{2g-1}$, and $\theta(G)\le 2 + \sqrt{2g-1}$ if $2g-1 = k^2$ for some integer $k$.
\end{theorem}
\pf Let $\Psi(G)$ be an embedding of $G$ in $S_g$.
 Without loss of generality we assume that $\Psi(G)$ is an orientable genus embedding.
 By Lemmas~\ref{lem-2-1} and \ref{lem-2-2}, we may assume that $\Psi(G)$ contains a spanning disk $D$.
Let $E_0$ be the set of edges contained in $D$.
Then,  $(V(G), E_0)$ is a planar graph.
We choose $(V(G), E_0)$ as one of our planar graphs in the edge partition for thickness.
The rest of the edge partition is part of the edge partition of $G^e$, the subgraph of $G$ consisting of essential edges of $\Psi$.

Let $d=2+\lfloor\sqrt{2g-1}\rfloor$, and let $G'$ be the graph obtained from $G^e$ by iteratively removing vertices of degree at most $d$ until no such vertex exists. Let $n'$ and $m'$ be the number of vertices and edges of $G'$.
The embedding of $G'$ is a subembedding of $\Psi(G)$ and no edge on $D$ is contained in $G'$.
 We need to add at least $n' +(n'-3)$ edges to make the subembedding  of $G'$ to an embedding of a spanning supergraph $G''$ of $G'$ with all faces being triangles.
 Here the $n'$ edges along the boundary of $D$ form a hamilton cycle $C$ of $G'$ and the $n'-3$ edges are  chords of $C$ embedded in $D$.
 Note that this is not necessarily a triangulation since $G''$ may have some pair of parallel edges in which one is an essential edge of $\Psi(G)$ and the other is the new edge of $G''$ on the cycle $C$.
 The purpose of adding edges/multiple edges is only for edge counting.
 Therefore, $m'+2n'-3\le 3n'-6+6g$ by Euler's formula and the fact that all faces are of size at least three.  Then, $\delta(G')\ge d+1\ge 2    +\sqrt{2g-1}$, and thus

$$(2+\sqrt{2g-1})n'\le 2m'\le 2n'-6+12g.$$

This inequality can be simplified as
$$n'\le 6 \sqrt{2g-1}.$$

Now, $\theta(G')\le \theta(K_{n'})\le \lfloor{n'+7\over 6}\rfloor<2+\sqrt{2g-1}$ by (\ref{thick-complete}). The integrity shows that $\theta(G')\le d$.

Let $Q_1, \ldots, Q_d$ be a decomposition of $E(G')$ into planar graphs.
Following Lemma 3.4, we decompose $E({G^e})\setminus E(G')$ into forests $F_1, \ldots, F_d$ such that, for each $i$, each component of $F_i$ has at most one  vertex in common with $G'$.
Now, $E_0\cup \{Q_i\cup F_i \;:\; 1\le i\le d\}$ is a decomposition of $E(G)$ into $d+1$ planar graphs,
and the theorem is true.

If $2g-1 = k^2$ for some integer $k$,
we let $d=1+ \sqrt{2g-1}$, then $d+1=2+ \sqrt{2g-1}$.
The statement $\theta (G) \le d+1 $ still holds, and thus $\theta (G) \le 2+ \sqrt{2g-1}$ as desired. \qed

\medskip

With similar arguments  to those above, we obtain an upper bound on the outerthickness of graphs.

\begin{theorem}\label{thm-3-5}
Let $G$ be a connected  graph embedded in $S_g$ ($g\ge 1$).
 Then, $\theta_o(G)\le 4 + \sqrt{3g-3/2}$,
 and $\theta_o(G)\le 3 + \sqrt{3g-3/2}$ if $3g-3/2 =h^2$ for some integer $h$.
\end{theorem}
\pf As in the proof of Theorem~\ref{thm-3-4}, we choose an embedding of $G$ such that all of its vertices are on a disk, and let $E_0$ be the set of edges on the disk. Now, $(V(G), E_0)$ is a planar graph, and thus has outerthickness 2 by Goncalves' conclusion [8].
 Let $G^e=G-E_0$.
 Then, in the embedding of $G^e$ induced by that of $G$, all  its vertices are on the boundary of a spanning disk,  and all its edges are essential.

Let $d=2+\lfloor\sqrt{3g-3/2}\rfloor$, and let $G'$ be the graph obtained from $G^e$ by iteratively removing vertices of degree at most $d$ until no such vertex exists. Then, $\delta(G')\ge d+1\ge 2+\sqrt{3g-3/2}$, and
$$(2+\sqrt{3g-3/2})|V(G')|\le 2|V(G')|-6+12g$$ by Euler's formula.

This inequality can be simplified to
$$|V(G')|\le 4 \sqrt{3g-3/2}.$$

Now, $\theta_o(G')\le \theta_o(K_{|V(G')|})\le \lceil{|V(G')|+2\over 4}\rceil<2+\sqrt{3g-3/2}$ by (\ref{outer-complete}), and the integrity shows that $\theta_o(G')\le d$.

Let $O_1, \ldots, O_d$ be a decomposition of $E(G')$ into outerplanar graphs, and let $F_1, \ldots, F_d$ be a decomposition of $E(G)\backslash E(G')$ into forests  such that, for each $i$, each component of $F_i$ has at most one  vertex in common with $G'$ (such a decomposition does exists by Lemma~\ref{lem-3-1}).
Let $O, O'$ be a decomposition of $E_0$ into two outerplanar graphs. Now, $\{O, O'\}\cup \{O_i\cup F_i \;:\; 1\le i\le d\}$ is a decomposition of $E(G)$ into $d+2$ outerplanar graphs.

If $3g-3/2 =h^2$ for some integer $h$, we let $d= 1+ \sqrt{3g-3/2}$. Then $d+1= 2+ \sqrt{3g-3/2}$,
and $\theta _o (G) \le d+2$ still holds.
This implies that $\theta_o(G)\le 3 + \sqrt{3g-3/2}$ if  $3g-3/2 =h^2$ for some integer $h$. \qed

\medskip

Using similar arguments as above, we have the following two theorems for nonorientable surfaces.

\begin{theorem}\label{thm-3-7}
Let $G$ be a connected  graph embedded in $N_k$ ($k\ge 1$).
  Then, $\theta(G)\le 3+ \sqrt{k-1}$,
  and $\theta (G)\le 2 + \sqrt{k-1}$ if $k-1 =h^2$ for some integer $h$.
\end{theorem}

\begin{theorem}\label{thm-3-8}
Let $G$ be a connected  graph embedded in $N_k$ ($k \ge 1$).
 Then, $\displaystyle \theta_o(G)\le 4 + \sqrt{ {3 \over 2} (k-1) }$,
 and $\theta_o(G)\le 3 + \sqrt{ {3 \over 2} (k-1) }$ if $ {3 \over 2} (k-1)=h^2$ for some integer $h$.
\end{theorem}

\medskip

\section{Outerthickness for toroidal graphs}

In previous section we provide some upper bounds of thickness and outerthickness for graphs in terms of their genus. Here, we will find the exact value for $\theta_o (S_1)$, the outerthickness of the torus.

Asano \cite{As87} proved  that $\theta (S_1) =2$.
Therefore $\theta_o (S_1) \le 4$ as each planar graph is partitioned into at most two outerplanar graphs \cite{Go05}. We will show in fact

\begin{theorem}\label{thm-4-1}
   $\theta_o(S_1) =3$.
\end{theorem}

We only need to show  $\theta_o(S_1) \le 3$ since  $\theta_o(K_7) =3$ and $K_7$ is toroidal.
We first prove a result for nonhomotopic loop systems for the torus.
Recall that $$\rho(\Sigma)= \max \{ |{\cal L}|: {\cal L} \mbox{ is a maximal nonhomotopic loop system of} \Sigma  \}.$$

\begin{lemma}\label{lem-4-2}
$\rho(S_1) =3$.
\end{lemma}
\pf Let ${\cal L}$ be a maximal nonhomotopic loop system of the torus with base point $x$.
Two noncontractible loops in the torus are homotopic if and only if they are homotopically disjoint and bound a cylinder.
All loops in ${\cal L}$ are nonhomotopic and they intersect only at the base point $x$.
This implies that all loops in ${\cal L}$ mutually intersect transversely.
Let $l_1, l_2, ..., l_t$ be loops contained in ${\cal L}$.
The fundamental group of torus is $Z \oplus Z$, generated by two elements which are not homotopic and intersect each other transversely.
Therefore,  up to homeomorphism, we may assume $l_1$ and $l_2$ are two generators of the torus.
View $l_1$ and $l_2$ as two edges of a graph with one vertex $x$ embedded in the torus.
The only face of this embedding has facial walk of  $x l_1 x l_2 x l_1 ^{-1} x l_2 ^{-1}$ which is a $4$-gon.
The vertex $x$ appears on this facial walk four times.
We may technically call these four appearance of $x$ as $x_1, x_2, x_3$ and $x_4$.
We can only connect $x_1 $ with $x_3$, or $x_2$ with $x_4$ to embed the third loop $l_3$.
After $l_3$ is embedded, there is no other way to embed any nonhomotopic loop inside the face $f$.
Therefore Lemma 4.2 is true.  \qed

\medskip

\begin{figure}[htp]
  \begin{center}
    \includegraphics[width=8.5cm]{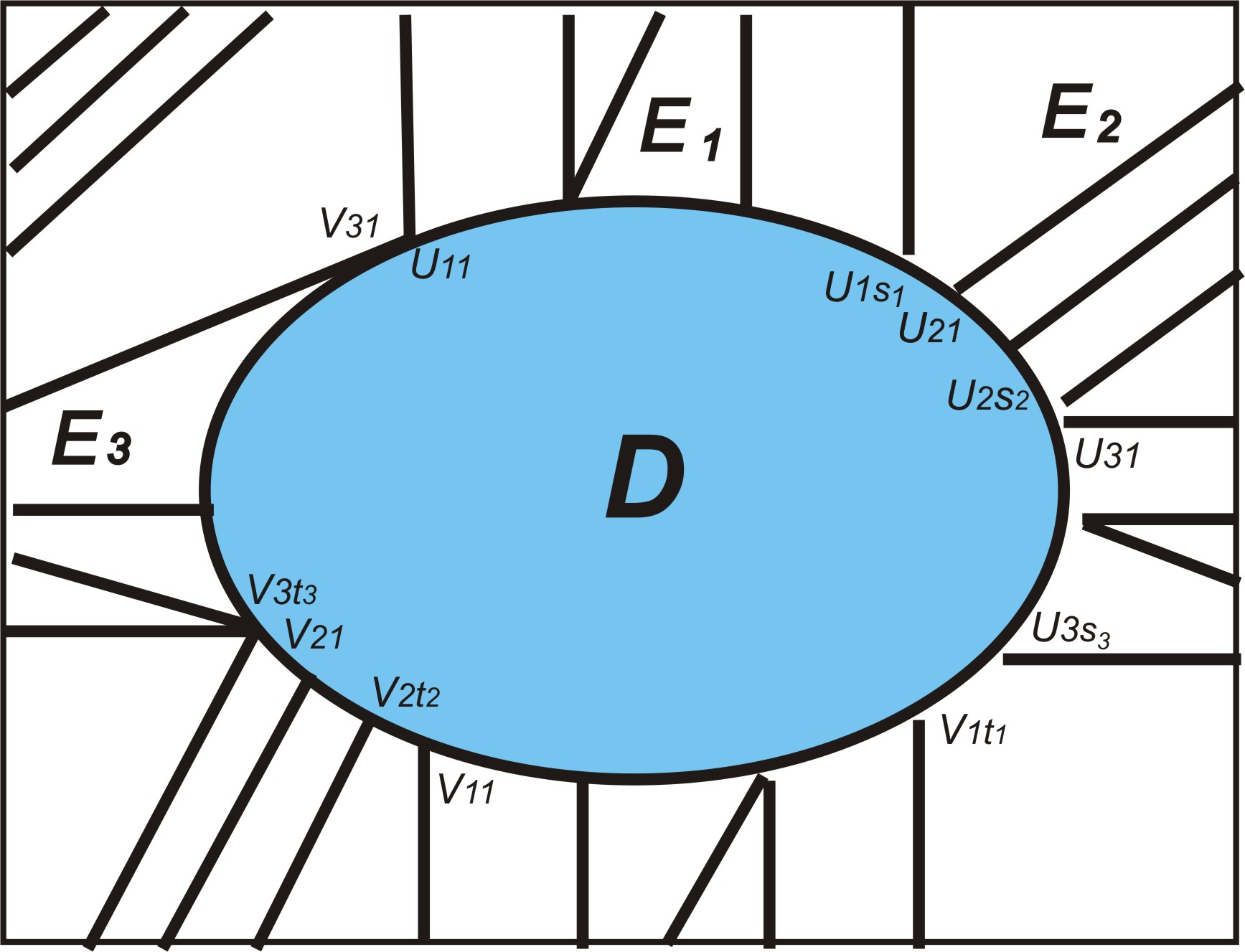} %
    \caption{Three homotopy classes of noncontractable loops of the torus.}
    \label{torus}
  \end{center}
\end{figure}

\noindent{\bf Proof of Theorem~\ref{thm-4-1}}.
Let $G$ be a graph embedded in the torus with embedding $\Psi (G)$.
By Lemmas~\ref{lem-2-1} and \ref{lem-2-2} we may assume that  $\Psi(G)$ contains a spanning disk $D$ which contains all vertices of $G$ and is bounded by a cycle $C$ of $G$.
All edges not contained in $D$ are essential edges.
If we contract $D$ to a point, then each essential edge is a noncontractable loop of the torus.
By Lemma~\ref{lem-4-2}  all essential edges are partitioned into at most three homotopy classes.
Certainly it is   easier to find the thickness and outerthickness for toroidal embeddings with one or two homotopy classes of essential edges than
with three homotopy classes.
Therefore we may assume that all three classes are not empty and they are embedded in the torus as shown in Figure~\ref{torus}.


\bigskip

\medskip

In order to partition the edge set of $G$ into three outerplanar graphs,
our strategy is to choose a spanning supergraph of $D$ ($D$ with some essential edges) as a planar graph and then partition this planar graph into
two outerplanar graphs, and use the remaining edges (all contained in $G^e$) as the third outerplanar subgraph.

Let  ${\cal E}_1, {\cal E}_2$ and ${\cal E}_3 $ be these three sets of essential edges with
 ${\cal E}_i =\{u_{ij}v_{ik} \}$ for some $j\in \{1, 2, ..., s_i \} $ and  $k\in \{1, 2, ..., t_i \}$.
 Then $G^e = {\cal E}_1 \cup {\cal E}_2 \cup {\cal E}_3 $.
All endvertices of essential edges are contained in the cycle $C$ in the clockwise order:
  $u_{11}...u_{1s_1}u_{21}...u_{2s_2}u_{31}...u_{3s_3}  v_{1t_1}...v_{11} $ $ v_{2t_2} ...v_{21} v_{3t_3} ... v_{31}$.
It is possible that $u_{1s_1} =u_{21}, u_{2s_2} =u_{31}, u_{3s_3} =v_{1t_1}, v_{11} =v_{2t_2}, v_{21} =v_{3t_3}$ and/or $v_{31} =u_{11}$.

For $i=1,2,3,$ let $G_i$ be the graph consisting of edges in ${\cal E}_i$.
We study the structure of $G_1$.
Similar results  also apply to $G_2$ and $G_3$.
A vertex $y$ of $G_1$ is an \it internal vertex \rm if $u_{11} \ne u_{1s_1}$ and $y \in u_{11} C u_{1s_1}$
    or $v_{11} \ne v_{1t_1}$ and $y \in v_{1t_1} C v_{11}$.
Define \it internal vertices \rm for $G_{2}$ and $G_{3}$ similarly.

We make some reductions.
We may delete all the internal vertices of $G_1$   with degree one since  adding a pendant  edge to a graph
does not increase the outerplanar  thickness.
Such deletions can be performed successfully.
An internal vertex of degree more than two must be incident to some pendant edges.
So we may suppose that all internal vertices of $G_1$ are of degree two.

If $G_1$ contains some internal vertices  then all internal vertices form a path joining $u_{11} $ to $u_{1s_1}$, or $u_{11} $ to $v_{1t_1}$,
     or $v_{11}$ to $u_{1s_1}$, or $v_{11}$ to $v_{1t_1}$ exclusively.
In this case we have (i) $u_{11} \ne u_{1s_1}$ and $v_{11} \ne v_{1t_1}$,
   (ii) $G_1$ is a path, and (iii) $G_1 \cup D$ is homeomorphic to a cylinder.

If $G_1$ does not contain any internal vertices,
then $G_1$  contains two essential edges $u_{11} v_{11}$ and $u_{1s_1} v_{1t_1}$ for which $u_{11}$ and $u_{1s_1}$ may coincide and/or
$v_{11}$ and $v_{1t_1}$ may coincide.
If $u_{11} \ne u_{1s_1}$ and $v_{11} \ne v_{1t_1}$ then $G_1$ contains  two disjoint edges $u_{11} v_{11}$ and $u_{1s_1} v_{1t_1}$,
with a possible third edge $u_{11}v_{1t_1}$ or $v_{11}v_{1s_1}$,
and $G_1 \cup D$ is homeomorphic to a cylinder.
If $u_{11} = u_{1s_1}$ and $v_{11}= v_{1t_1}$, then $G_1$ consists of only one single essential edge $u_{11} v_{11}$.
If $u_{11} = u_{1s_1}$ and $v_{11} \ne  v_{1t_1}$ or $u_{11} \ne u_{1s_1}$ and $v_{11}= v_{1t_1}$,
then $G_1$ consists of two essential edges forming a path of length two, and $G_1 \cup D$ is homeomorphic to a degenerate cylinder.
If $G_1$ consists of  only two disjoint edges, we add a new edge $u_{11} u_{1s_1}$ to connect these two edges to make the new $G_1$ again a path. Therefore $G_1$ is just a path in any case.

If one of $G_1 \cup D, G_2 \cup D$ and $G_3 \cup D$, say $G_1 \cup D$, is homeomorphic to a cylinder (not a degenerate cylinder, i.e., $u_{11} \ne u_{1s_1}$ and $v_{11} \ne v_{1t_1}$),
then $u_{21} \ne v_{31} $ and $u_{3s_3} \ne v_{2t_2}$.
Hence $G_2$ and $G_3$ can possibly  have at most two common vertices which are $u_{2s_2}$ (= $u_{31}$) and $v_{21} $ (= $v_{3t_3}$).
By the structure of $G_2$ and $G_3$ (similar to $G_1$,   since each is a path),
$G_2 \cup G_3$ is the union of two paths that have up to two common vertices, hence is outerplanar.
Since $D \cup G_1$ is a planar graph, by Gonalves \cite{Go05}, $D\cup G_1$ can be partitioned into two outerplanar subgraphs $H_1$ and $H_2$.
Therefore $G$ can be partitioned into three outerplanar graphs $H_1, H_2$ and $G_2 \cup G_3$.

If none of $G_1 \cup D, G_2 \cup D$ and $G_3 \cup D$ is homeomorphic to a cylinder,
then  $G_i, i=1, 2, 3$, is either a single edge or a path of length 2.
Clearly $G_2 \cup G_3$ is an outerplanar graph, and by similar reasoning, $G$ can be partitioned into three outerplanar graphs. This completes the proof of Theorem~\ref{thm-4-1}. \qed

\section{Concluding remarks}


The thickness and outerthickness problems are ${\cal N}{\cal P}$-hard and therefore there are not many results on the exact value for graphs other than a few classes of graphs with high symmetry.
Su, Kanno and VanHeeswijk \cite{SKV} asked whether all projective planar graphs have outerthickness  $2$.
This is certainly an interesting question after Concalves confirmed Chartrand et al.'s conjecture.
While no counterexample has been found, it seems much more difficult to show that all projective planar graphs have outerthickness $2$ because all nonplanar and  projective planar graphs have thickness $2$.
Also there is even no statement to claim the thickness of the projective plane is $2$.
We are able to show that the thickness of the projective plane and the Klein bottle both are $2$ and the outerthickness of the Klein bottle is $3$.
We will include these results, mainly for graphs embedded in nonorientable surface, in another paper.

We note that both upper bounds for thickness and outerthickness by Dean and Hutchinson, and by Theorems~\ref{thm-3-4}-\ref{thm-3-8} are correct in the order
of $O(\sqrt{g})$ or $O(\sqrt{k})$.
This can be explained by the thickness and outerthickness of the complete graphs $K_n$.
We know that $\theta (K_n) = \lfloor (n+7)/6 \rfloor$ and $\theta_o (K_n) = \lceil (n+1)/4 \rceil$ by (\ref{thick-complete}) and (\ref{outer-complete}).
The upper bound provided by Theorem~\ref{thm-3-4} is
  $\displaystyle 3 + \sqrt{2\lceil {(n-3)(n-4) \over 12 } \rceil -1  }$ $\le {1\over \sqrt{6}} n +3$.

\end{document}